\documentclass[pamm,a4paper,fleqn]{w-art}
\usepackage{times,cite,w-thm}
\usepackage[T1]{fontenc}
\usepackage[utf8]{inputenc}
\theoremstyle{plain}

\theoremstyle{definition}

\numberwithin{equation}{section}
\usepackage{amssymb}
\usepackage{tikz}
\usepackage{graphicx}

\newcommand\C{{\mathbb C}}
\newcommand\N{{\mathbb N}}

\newcommand\R{{\mathbb R}}
\newcommand\T{{\mathbb T}}

\newcommand\Z{{\mathbb Z}}

\newcommand\eps{{\varepsilon}}
\newcommand\spec{{\rm spec}\,}  % with space
\newcommand\specn{{\rm spec}}   % no space
\newcommand\Specn{{\rm Spec}}   % no space
\newcommand\speps{{\rm spec}_\eps}
\newcommand\Speps{{\rm Spec}_\eps}
\newcommand\spess{{\rm spec}_{\rm ess}\,}
\newcommand\diam{{\rm diam}}
\newcommand\supp{{\rm supp}\,}

\begin{document}
%% \def\leftmark{Session title}
%%
%%    The information for the title page will be placed between
%%    \begin{document} and \maketitle. The order of most entries
%%    is determined by the class file and can not be changed by
%%    rearranging them. The maketitle command follows after the
%%    abstract.
%%
%%    Most of the following commands will be completed by the publisher.
%%
%%    \renewcommand{\copyrightyear}{2016}
%%    \DOIsuffix{pamm.20161zzzz}
%%    \Volume{16} 
%%    \Year{2016} 
%%    \pagespan{1}{}
%%
%%    The short title is optional:

\TitleLanguage[EN]
\title{Convergent spectral inclusion sets for banded matrices}

%% Please do not enter footnotes or \inst{}-notes into the optional
%% argument of the author command. 

%% Please delete not needed author entries.
%% Information for the first author.
\author{\firstname{Simon N.}  \lastname{Chandler-Wilde}\inst{1,}%
  \footnote{s.n.chandler-wilde@reading.ac.uk}}

\address[\inst{1}]{\CountryCode[GB] Department of Mathematics and Statistics, University of Reading, RG6 6AX, UK.}
%%
%%    Information for the second author
\author{\firstname{Ratchanikorn}  \lastname{Chonchaiya}\inst{2,}%
  \footnote{hengmath@hotmail.com}}

\address[\inst{2}]{\CountryCode[TH] Faculty of Science, King Mongkut's University of Technology Thonburi, Bangkok 10140, Thailand.}
%%
%%    Information for the third author
\author{\firstname{Marko} \lastname{Lindner}\inst{3,}%
  \footnote{Corresponding author: \ElectronicMail{lindner@tuhh.de}}}

\address[\inst{3}]{\CountryCode[DE] Institut Mathematik, TU Hamburg, 21073 Hamburg, Germany.}
%%
%%    \dedicatory{This is a dedicatory.}
%%
%%    Abstract is required.
\AbstractLanguage[EN]
\begin{abstract}
{\bf Abstract.} 
We obtain sequences of inclusion sets for the spectrum, essential spectrum, and pseudospectrum of banded, in general non-normal, matrices of finite or infinite size. %In the infinite case (bi- or semi-infinite), the matrix acts as a bounded linear operator on the corresponding $\ell^2$ space, and we moreover bound and approximate its essential spectrum. 
Each inclusion set is the union of the pseudospectra of certain submatrices of a chosen size $n$. Via the choice of $n$, one can balance accuracy of approximation against computational cost, and we show, in the case of infinite matrices, convergence as $n\to\infty$ of the respective inclusion set to the corresponding spectral set.
\end{abstract}
%% maketitle must follow the abstract.
\maketitle                   % Produces the title.

%%%%%%%%%%%%%%%%%%%%%%%%%%%%%%
\section{Introduction}
In many finite difference schemes or in physical or social models, where interaction between objects is direct in a finite radius only (and is of course indirect on a global level), the corresponding matrix or operator is {\sl banded}, also called {\sl of finite dispersion}, meaning that the matrix is supported on finitely many diagonals only. 
%\[
%%\text{\Large $A$}  
%A\ =\ 
%\begin{pmatrix}~
%\begin{tikzpicture}[scale=0.53]
%\foreach \k in {0,0.5,...,4.5}
%  \fill[gray!50] (\k,6-\k) rectangle ++(1.5,-1.5);
%\draw[gray!30,very thin,step=0.5cm] (0,0) grid (6,6);
%\end{tikzpicture}
%~\end{pmatrix}
%\]
In the case of finite matrices this is of course a tautology; in that contect one assumes that the {\sl bandwidth} is not only finite but small compared to the matrix size, where the bandwidth of a matrix $A$ is the distance from the main diagonal in which nonzeros can occur. (Precisely: it is the largest $|i-j|$ over all matrix positions $(i,j)$ with $A_{i,j}\ne 0$.)
So this is our setting: finite, semi-infinite or bi-infinite banded matrices.

We equip the underlying vector space with the Euclidian norm, so our operators act on an $\ell^2$ space over $\{1,\ldots,N\}$ or $\N=\{1,2,\ldots\}$ or $\Z=\{\ldots,-1,0,1,\ldots\}$.
In the two latter cases (semi- and bi-infinite matrices), we assume each diagonal to be a bounded sequence, whence the matrix acts as a bounded linear operator, again denoted $A$, on the corresponding $\ell^2$ space.

The exact computation of the spectrum by analytical means is in general impossible (by Abel-Ruffini) if the size of the matrix is larger than four. So one is forced to resort to approximations. But for non-normal matrices and operators, also the approximation of the spectrum is extremely delicate and unreliable, hence one often substitutes for the spectrum, $\spec A$, the {\sl pseudospectrum},
\[
\speps A\ :=\ \{\lambda\in\C:\|(A-\lambda I)^{-1}\|> 1/\eps\},\qquad \eps>0,
\]
that is much more stable to approximate, and then sends $\eps\to 0$. Note that, by agreeing to say $\|B^{-1}\|=\infty$ if $B$ is not invertible, one has $\spec A\subset\speps A$ for all $\eps>0$.
For an impressive account of pseudospectra and their applications, see the monograph \cite{TrefEmb}.

Our aim in this paper is to derive inclusion sets for 
%show %(and to some extent explain, see \cite{SCW-ML:Inclusion1,SCW-ML:Inclusion2} for the full details) 
%subsets, supersets and approximations of 
$\spec A,\ \speps A$ as well as the essential spectrum, $\spess A$, in terms of unions of pseudospectra of moderately sized (but many) finite submatrices of $A$ of column dimension $n$. Moreover, if the matrix is infinite, we prove convergence, as $n\to\infty$, of the respective inclusion set to each of $\spec A$, $\speps A$, or $\spess A$.

\section{Approximating the lower norm on $\ell^2(\Z)$}
Our arguments are, perhaps surprisingly, tailor-made for the case  of bi-infinite vectors and matrices on them. In fact, instead of $\ell^2(\Z)$, everything also works for $\ell^2(G)$ with a discrete group $G$, e.g. $G=\Z^d$, subject to Yu's so-called Property A \cite{HagLiSei,Yu}. Only later, in Section \ref{sec:semi+finite}, we manage to work around the group structure and to transfer results to $\ell^2(\N)$ and $\ell^2(\{1,\ldots,N\})$, hence: to semi-infinite and finite matrices.

As some sort of antagonist of the operator norm, $\|A\|=\sup\{\|Ax\|:\|x\|=1\}$, we look at
the so-called {\sl lower norm}\footnote{It is not a norm! Our terminology is that of  \cite{RaRoSiBook,MLBook}.}
\[
\nu(A)\ :=\ \inf\{\|Ax\|:\|x\|=1\}
\]
of a banded and bounded operator on $\ell^2(\Z)$. Fixing $n\in\N$ and limiting the selection of unit vectors $x$ to those with a finite support of diameter less than $n$, further limits how small $\|Ax\|$ can get. Precisely,
\begin{equation} \label{eq:nundef}
\nu_n(A)\ :=\ \inf\{\|Ax\|:\|x\|=1,\ \diam(\supp x)<n\},\qquad n\in\N,
\end{equation}
is typically larger than $\nu(A)$ -- but (and this is remarkable) only larger by at most the amount of
a certain $\eps_n\sim 1/n$ that we will quantify precisely below. Let us first write this important fact down:
\begin{equation} \label{eq:nun}
\nu_n(A)-\eps_n\ \le\ \nu(A)\ \le\ \nu_n(A).
\end{equation}
This observation can be traced back to \cite{HengPhD, BigQuest} %,SCW-ML:Inclusion1,SCW-ML:Inclusion2}
and, for Schrödinger operators, even to \cite{Elliot,AvronMoucheSimon}. 
Extensive use has been made of \eqref{eq:nun}, e.g., in \cite{TriRand,subwords}.
%
%The quantity $\eps_n$ depends on $n$, $\|A\|$ and the bandwidth of $A$. In particular, if we replace $A$ by an operator $B$ on $\ell^2(\Z)$ with norm and bandwidth not exceeding those of $A$ then
%\begin{equation} \label{eq:nun}
%\nu_n(A)-\eps_n\ \le\ \nu(A)\ \le\ \nu_n(A)
%\end{equation}
%also applies to $B$ in place of $A$ -- showing that this localizability of $\nu(\cdot)$ is rather a property of $\ell^2$ than of an operator on it.
%
The statement $\diam(\supp x)<n$ in \eqref{eq:nundef} translates to $\supp x\subseteq k+\{1,\ldots,n\}$ for some $k\in\Z$. Hence,
\begin{equation} \label{eq:nun_inf}
\nu_n(A)\ =\ \inf\{\nu(A|_{\ell^2(k+\{1,\ldots,n\})}):k\in\Z\}.
\end{equation}
Let us write $P_{n,k}$ for the operator of multiplication by the characteristic function of $k+\{1,\ldots,n\}$ and agree on writing
\[
A|_{\ell^2(k+\{1,\ldots,n\})}\ =:\ AP_{n,k}\ :\ \ell^2(k+\{1,\ldots,n\})\ \to\ \ell^2(\Z),\qquad n\in\N,\ k\in\Z.
\]
In matrix language, $AP_{n,k}$ corresponds to the matrix formed by columns number $k+1$ to $k+n$ of $A$. By the band structure of $A$, that submatrix is supported in finitely many rows only, even reducing it to a finite $m\times n$ matrix, where $m$ equals $n$ plus two times the bandwidth of $A$. Then $\nu(AP_{n,k})$, as in \eqref{eq:nun_inf}, is the smallest singular value of this $m\times n$ matrix, making this a standard computation.

\section{$\eps_n$ and the reduction to tridiagonal form} \label{sec:tridiag}
Our analysis of $\eps_n$ is particularly optimized in the case of tridiagonal matrices, that is when $A$ has bandwidth one, so that it is only supported on the main diagonal and its two adjacent diagonals. Let $\alpha,\beta,\gamma\in\ell^\infty(\Z)$ denote, in this order, the sub-, main- and superdiagonal of $A$, with entries $\alpha_i=A_{i+1,i},\ \beta_i=A_{i,i}$ and $\gamma_i=A_{i-1,i}$ with $i\in\Z$. In that case (see \cite{HengPhD,SCW-ML:Inclusion1,SCW-ML:Inclusion2}),
\begin{equation} \label{eq:epsn}
\eps_n\ =\ 2(\|\alpha\|_\infty+\|\gamma\|_\infty)\sin\frac \pi{2(n+1)}\ <\ (\|\alpha\|_\infty+\|\gamma\|_\infty)\frac \pi{n+1}\ \sim\ \frac 1n.
\end{equation}
Although \eqref{eq:nun} holds with this choice of $\eps_n$ for the very general setting of all tridiagonal matrices, formula \eqref{eq:epsn} turns out to be best possible in some nontrivial examples such as the shift operator \cite{SCW-ML:Inclusion1}.

To profit from these well-tuned parameters also in the case of larger bandwidths, note that \eqref{eq:nun} and \eqref{eq:epsn} even work in the block case, that is when the entries in the $\ell^2$ vectors are themselves elements of some Banach space $X$ and the matrix entries of $A$ are operators on $X$. So the trick with a band matrix $B$ with a larger bandwidth $b$ is to interpret $B$ as block-tridiagonal with blocks of size $b+1$:
\[
%\text{\Large $A$}  
B\ =\ 
\begin{pmatrix}~
\begin{tikzpicture}[scale=0.53]
\foreach \k in {0,0.5,...,4.5}
  \fill[gray!50] (\k,6-\k) rectangle ++(1.5,-1.5);
\draw[gray!30,very thin,step=0.5cm] (0,0) grid (6,6);
\end{tikzpicture}
~\end{pmatrix}
\ \cong\ 
\begin{pmatrix}~
\begin{tikzpicture}[scale=0.53]
\foreach \k in {0,0.5,...,4.5}
  \fill[gray!50] (\k,6-\k) rectangle ++(1.5,-1.5);
\draw[gray!30,very thin,step=0.5cm] (0,0) grid (6,6);
\draw[black,line width=0.5mm,step=1.5cm] (0,0) grid (6,6);
\foreach \i in {1,...,4}
  {
  \foreach \j in {1,...,4}
     \node at ({1.5*\j-0.75},{6.75-1.5*\i}) {$A_{{\i},{\j}}$};
  }
\end{tikzpicture}
~\end{pmatrix}
\ =\ A.
%\ =\ 
%\begin{pmatrix}~
%\begin{tikzpicture}[scale=0.53]
%\draw[gray!10,very thin,step=0.5cm] (0,0) grid (6,6);
%\draw[gray,line width=0.3mm,step=1.5cm] (0,0) grid (6,6);
%\foreach \k [evaluate=\k as \l using int(80-20*\k)] in {1,2,3}
%  {
%    \fill[gray!\l] ({1.5*(\k-1)},{7.5-1.5*\k}) rectangle ++(3,-3);
%    \draw[black,line width=0.7mm] ({1.5*(\k-1)},{7.5-1.5*\k}) rectangle ++(3,-3); 
%    \node at ({1.5*\k-0.5},{6.5-1.5*\k}) {\Large $S_{\k}$};
%  }
%\end{tikzpicture}
%~\end{pmatrix}
\]
Here, a matrix $B$ with bandwidth $b=2$ is identified with a block-tridiagonal matrix $A$ with $3\times 3$ blocks, noting that $3=b+1$.

Since the blocks of $A$ can be operators on a Banach space $X$, one can even study $\spec B$ and $\speps B$ by our techniques for bounded operators $B$ on $L^2(\R)\cong\ell^2(\Z,X)$, where $X=L^2([0,1])$, e.g.~for integral operators $B$ with a banded kernel $k(\cdot,\cdot)$.

\section{The role of the lower norm in spectral computations}
If $A$ sends a unit vector $x$ to a vector $Ax$ with norm $\frac 14$ then, clearly, $A^{-1}$, bringing $Ax$ back to $x$, has to have at least norm $4$. 
The lower norm, $\nu(A)$, is pushing this observation to the extreme. By minimizing $\frac{\|Ax\|}{\|x\|}$, it minimizes $\frac{\|y\|}{\|A^{-1}y\|}$ and hence computes the reciprocal of $\|A^{-1}\|$ -- with one possible exception: non-invertibility of $A$ due to $\nu(A)=0$ or $\nu(A^*)=0$.
%The lower norm is looking for the smallest (infimum) case in the forward scenario, hence directly hinting at $\|A^{-1}\|$ and from there at the resolvent norm, at spectra and pseudospectra. 
Properly: Since $\nu(A)>0$ iff $A$ is injective and has a closed range (a standard result), $A$ is invertible iff both $\nu(A)$ and $\nu(A^*)$ are nonzero\footnote{By $A^*$ we denote the Banach space adjoint of $A$. In particular, $(\lambda I)^*=\lambda I$, not $\overline\lambda I$.}. Keeping this symmetry of $A$ and $A^*$ in mind,
\[
1/\|A^{-1}\|\ =\ \min(\nu(A),\nu(A^*))\ =:\ \mu(A),
\]
(see, e.g., \cite{HagLiSei}), where, again, $\|A^{-1}\|=\infty$ signals non-invertibility and where $1/\infty:=0$. From here it is just a small step to
\[
\spec A\ =\ \{\lambda\in\C:\|(A-\lambda I)^{-1}\|=\infty\}\ =\ \{\lambda\in\C:\mu(A-\lambda I)=0\}
\]
\pagebreak
and
\begin{equation} \label{eq:speps}
\speps A\ =\ \{\lambda\in\C:\|(A-\lambda I)^{-1}\|> 1/\eps\}\ =\ \{\lambda\in\C:\mu(A-\lambda I)<\eps\},\qquad\eps>0.
\end{equation}
Being able to approximate $\nu(A)$, up to $\eps_n\sim \frac 1n$, by $\nu_n(A)$, enables us to approximate $\spec A$ and $\speps A$, with a controllable error, by sets built on $\nu_n(A-\lambda I)$ and $\nu_n(A^*-\lambda I)$.

\section{Approximating the pseudospectrum in the bi-infinite case}
Applying \eqref{eq:nun} to $A-\lambda I$ and $(A-\lambda I)^*=A^*-\lambda I$ in place of $A$, we see that
\[
\mu_n(A-\lambda I)<\eps\quad \implies\quad \mu(A-\lambda I)<\eps\quad\implies\quad\mu_n(A-\lambda I)<\eps+\eps_n,
\]
where $\mu_n(B):=\min(\nu_n(B),\nu_n(B^*))$, noting that $\eps_n$ is independent of $\lambda\in\C$, by \eqref{eq:epsn}.
Combining this with \eqref{eq:nun_inf} and \eqref{eq:speps}, we conclude (cf. \cite[Thm. 4.3 \& Cor. 4.4]{HengPhD}):

\begin{proposition} {\bf (bi-infinite case)}\label{prop:2sided}
For bounded band operators $A$ on $\ell^2(\Z)$ and corresponding $\eps_n$ from \eqref{eq:epsn}\footnote{Note that $\eps_n$, if using \eqref{eq:epsn}, has to be computed for the block-tridiagonal representation of $A$, see Section \ref{sec:tridiag}.}, one has
\begin{equation} \label{eq:sandwich}
\bigcup_{k\in\Z}\speps (AP_{n,k},A^*P_{n,k})
\quad\subseteq\quad \speps A\quad\subseteq\quad
\bigcup_{k\in\Z}\specn_{\eps+\eps_n} (AP_{n,k},A^*P_{n,k}),\qquad \eps>0,\ n\in\N,
\end{equation}
where we abbreviate $\speps (A,B)\ :=\ \speps A\ \cup\ \speps B$.
\end{proposition}

By iterated application of \eqref{eq:sandwich}, one can extend \eqref{eq:sandwich} to the left and right as follows:
\[
\specn_{\eps-\eps_n}A\ \subseteq\ 
\bigcup_{k\in\Z}\speps (AP_{n,k},A^*P_{n,k})
\ \subseteq\ \speps A\ \subseteq\ 
\bigcup_{k\in\Z}\specn_{\eps+\eps_n} (AP_{n,k},A^*P_{n,k})
\ \subseteq\ \specn_{\eps+\eps_n}A.
\]
And now, sending $n\to\infty$, we have $\eps_n\to 0$, by \eqref{eq:epsn}, and then Hausdorff-convergence $\specn_{\eps+\eps_n}A\ \to\ \speps A$ as well as $\specn_{\eps-\eps_n}A\ \to\ \speps A$, see e.g. \cite{Shargo}\footnote{In the case of a Banach space valued $\ell^2$, that Banach space should be finite-dimensional or subject to the conditions in Theorem 2.5 of \cite{Shargo}.}.
We conclude  (cf. \cite[Sec. 4.3]{HengPhD}):

\begin{proposition}
The subsets and supersets of $\speps A$ in \eqref{eq:sandwich} both Hausdorff-converge to $\speps A$ as $n\to\infty$.
\end{proposition}

\section{Approximating the pseudospectra of semi-infinite and finite matrices} \label{sec:semi+finite}
Now take a bounded and banded operator $A$ on $\ell^2(\N)$.
In\cite{SCW-ML:Inclusion2} we show how to reduce this case (via embedding $A$ into a bi-infinite matrix plus some further arguments) to the bi-infinite result:
\begin{proposition} {\bf (semi-infinite case)}\label{prop:1sided}
For bounded band operators $A$ on $\ell^2(\N)$ and corresponding $\eps_n$ from \eqref{eq:epsn}, one has
%\begin{equation} \label{eq:sandwich1}
\[
\bigcup_{k\in\N_0}\speps (AP_{n,k},A^*P_{n,k})
\quad\subseteq\quad \speps A\quad\subseteq\quad
\bigcup_{k\in\N_0}\specn_{\eps+\eps_{n}} (AP_{n,k},A^*P_{n,k}),\qquad \eps>0,\ \ n\in\N,
\]
%\end{equation}
where again $\speps (A,B) := \speps A \cup \speps B$.
Also here the sub- and supersets Hausdorff-converge to $\speps A$ as $n\to\infty$.
\end{proposition}

The technique that helps to deal with one endpoint on the axis can essentially be repeated for a second endpoint:
\begin{proposition} {\bf (finite case)}\label{prop:finite}
For finite band matrices $A$ on $\ell^2(\{1,\ldots,N\})$ with some $N\in\N$, one has
%\begin{equation} \label{eq:sandwich_fin}
\[
\bigcup_{k=0}^{N-n}\speps (AP_{n,k},A^*P_{n,k})
\quad\subseteq\quad \speps A\quad\subseteq\quad
\bigcup_{k=0}^{N-n}\specn_{\eps+\eps_{n}} (AP_{n,k},A^*P_{n,k}),\qquad \eps>0,\ \ 1\le n\le N,
\]
%\end{equation}
where again $\speps (A,B) := \speps A \cup \speps B$.
\end{proposition}
This time, of course, there is no way of sending $n\to\infty$, hence no Hausdorff-convergence result.

\section{Approximating spectra}
So far we have convergent subsets and supersets of $\speps A$ for $\eps>0$. The spectrum, $\spec A$, can now be Hausdorff-approximated via sending $\eps\to 0$. However, there is a more direct approach:
Introducing closed-set versions of pseudospectra,
\[
\Speps A\ :=\ \{\lambda\in\C:\|(A-\lambda I)^{-1}\|\ge 1/\eps\}\ =\ \{\lambda\in\C:\mu(A-\lambda I)\le\eps\},\qquad\eps\ge 0,
\]
we can prove identical copies of Propositions \ref{prop:2sided}, \ref{prop:1sided} and \ref{prop:finite}  with upper-case (i.e.~closed) instead of lower-case (i.e.~open) pseudospectra everywhere -- and including the case $\eps=0$, see \cite{SCW-ML:Inclusion2}. The latter brings convergent supersets  for $\spec A=\Specn_0 A$ right away, without the need for a further limit $\eps\to 0$.
Here is the new formula for the bi-infinite case, evaluated for $\eps=0$.
\begin{equation} \label{eq:sandwich0}
\bigcup_{k\in\Z}\spec (AP_{n,k},A^*P_{n,k})
\quad\subseteq\quad \spec A\quad\subseteq\quad
\bigcup_{k\in\Z}\Specn_{\eps_n} (AP_{n,k},A^*P_{n,k}),\qquad n\in\N.
\end{equation}

\section{Examples}
For three selected operator examples, we show the Hausdorff-convergent (as $n\to\infty$) superset bounds on $\spec A$ from \eqref{eq:sandwich0}. All three operators are given by tridiagonal bi-infinite matrices. Moreover, all three matrices are periodic, so that we can analytically compute the spectrum by Floquet-Bloch; that is, treating the $3$-periodic matrix as a $3\times 3$-block convolution on $\ell^2(\Z,\C^3)$ and turning that, via the corresponding block-valued Fourier-transform, into a $3\times 3$-block multiplication on $L^2(\T,\C^3)$, whose spectrum is obvious, see, e.g., Theorem 4.4.9 in \cite{Davies}.
For comparison, the exact spectrum is superimposed in each example as a red curve in the last column. 
\medskip

{\bf a) } We start with the right shift, where the subdiagonal is $\alpha=(\dots, 1,1,1,\dots)$ and the main and superdigonal are $\beta=\gamma=(\dots,0,0,0,\dots)$. The spectrum is the unit circle, and here are our supersets for $n\in\{4,8,16\}$:

\noindent
\begin{tabular}{ccc}
\includegraphics[width=55mm]{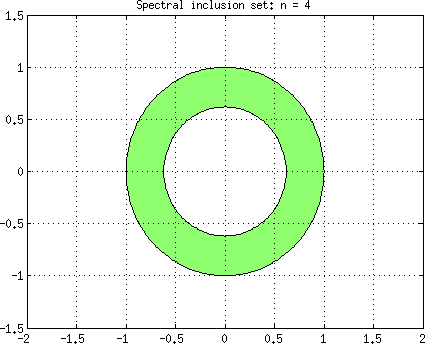}&
\includegraphics[width=55mm]{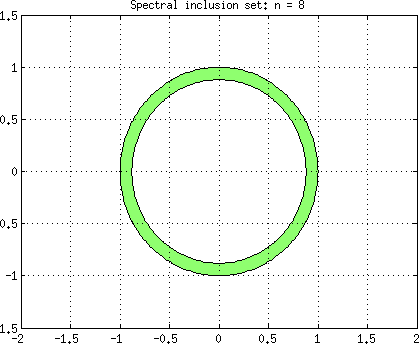}&
\includegraphics[width=55mm]{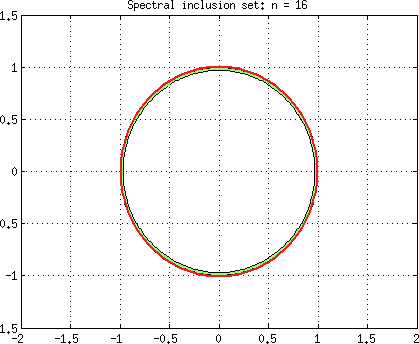}
\end{tabular}
\medskip

{\bf b) } Our next example is 3-periodic with subdiagonal $\alpha=(\dots,\mathbf 0,0,0,\dots)$,
main diagonal $\beta=(\dots,\mathbf{-\frac 32},1,1,\dots)$ and superdiagonal $\gamma=(\dots,\mathbf 1,2,1,\dots)$, where $\alpha_0,\ \beta_0$ and $\gamma_0$ are highlighted in boldface. The spectrum consists of two disjoint loops, and we depict our supersets for $n\in\{32, 64, 128\}$:

\noindent
\begin{tabular}{ccc}
\includegraphics[width=55mm]{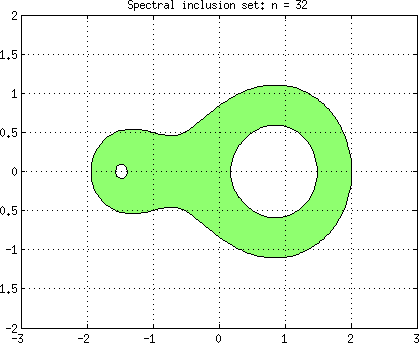}&
\includegraphics[width=55mm]{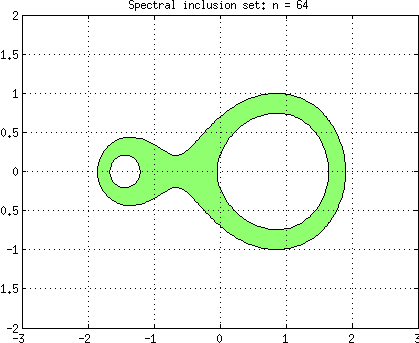}&
\includegraphics[width=55mm]{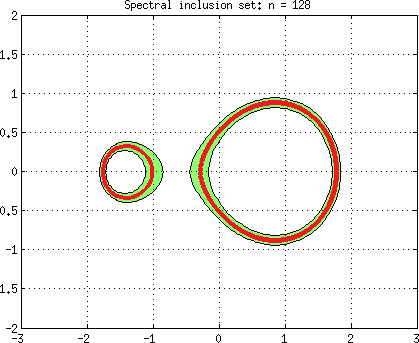}
\end{tabular}
\medskip

{\bf c) } Our third example is also 3-periodic with subdiagonal $\alpha=(\dots,\mathbf 0,0,0,\dots)$, main diagonal $\beta=(\dots,\mathbf{-\frac 12},1,1,\dots)$ and superdiagonal $\gamma=(\dots,\mathbf 1,2,1,\dots)$, where $\alpha_0,\ \beta_0$ and $\gamma_0$ are highlighted in boldface. The spectrum consists of one loop, and we depict our supersets for $n\in\{8,16,32\}$:

\noindent
\begin{tabular}{ccc}
\includegraphics[width=55mm]{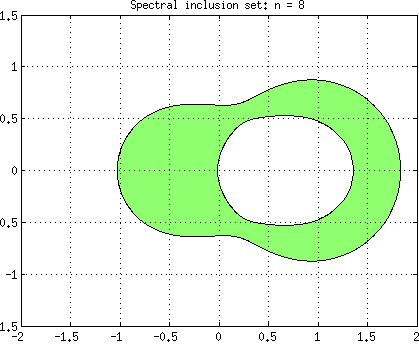}&
\includegraphics[width=55mm]{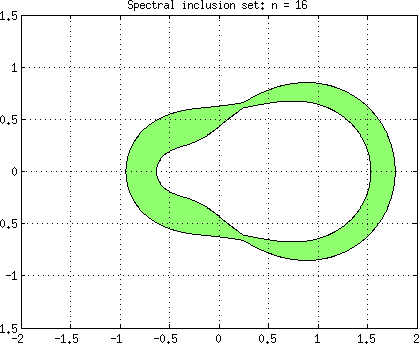}&
\includegraphics[width=55mm]{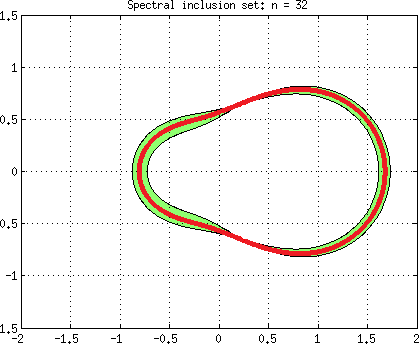}
\end{tabular}
\medskip

Another effect of the $3$-periodicity of the diagonals in $A$ is that there are only three distinct submatrices $AP_{n,k}$ and $A^*P_{n,k}$ each for $k\in\Z$.
In fact, for many operator classes, the infinite unions in \eqref{eq:sandwich}, \eqref{eq:sandwich0} and so on, reduce to finite unions. %via a-priori knowledge of the model. 
E.g., for a $\{0,1\}$-valued aperiodic diagonal \cite{LenzPhD}, there are only $n+1$ different subwords of length $n$, and for a $\{0,1\}$-valued random diagonal, there are $2^n$ (again, finitely many) different subwords of length $n$. Also, for non-discrete diagonal alphabets, the infinite union can be reduced to a finite one via compactness arguments, see our discussion in \cite{SCW-ML:Inclusion1}.

\section{Approximating essential spectra}
In the case where $A$ is an infinite matrix there is large interest also in the approximation of the {\sl essential spectrum}, $\spess A$, which is the spectrum in the Calkin algebra,
i.e.~ the set of all $\lambda\in\C$ where $A-\lambda I$ is not a Fredholm operator, i.e.~is not invertible modulo compact operators. 

Our results in this section apply when each $A_{i,j}\in\C$, but also when each $A_{i,j}$ is a bounded linear operator on a Banach space $X$, as long as $X$ is finite-dimensional or the operators $\{A_{i,j}\}$ are collectively compact in the sense of \cite{CollComp1,CollComp2}.
\medskip

As for the spectrum (see Section \ref{sec:tridiag}) it is enough to consider the case when $A$ is tridiagonal. 
%
%The question about the essential spectrum is hence only interesting for infinite band matrices.
The bi-infinite case is easily reduced to the semi-infinite case: 
Indeed, modulo compact operators, 
\[
A \ \cong\  
\left(\begin{array}{ccc|c|ccc}
\ddots&\ddots&&&&\\
\ddots&A_{-2,-2}&A_{-2,-1}&&\\
&A_{-1,-2}&A_{-1,-1}&&\\\hline
&&&0\\\hline
&&&&A_{1,1}&A_{1,2}\\
&&&&A_{2,1}&A_{2,2}&\ddots\\
&&&&&\ddots&\ddots
\end{array}\right)
\ =:\ 
\left(\begin{array}{c|c|c}
A_{-}&&\\\hline
&0&\\\hline
&&A_{+}
\end{array}\right),
\]
so that
\[
\spess A\ =\ \spess A_{-}\ \cup\ \spess A_{+}.
\]

It remains to look at semi-infinite banded matrices $A$.
Modulo compact operators, for every $m\in\N$,
\[
A\ =\ 
\begin{pmatrix}
A_{11}&A_{12}\\
A_{21}&A_{22}&\ddots\\
&\ddots&\ddots
\end{pmatrix}
\ \cong\ 
\left(\begin{array}{ccc|ccc}
0&&\\
&\ddots&\\
&&0\\\hline
&&&A_{m+1,m+1}&A_{m+1,m+2}\\
&&&A_{m+2,m+1}&A_{m+2,m+2}&\ddots\\
&&&&\ddots&\ddots
\end{array}\right),
\]
so that, with
\[
\begin{pmatrix}
A_{m+1,m+1}&A_{m+1,m+2}\\
A_{m+2,m+1}&A_{m+2,m+2}&\ddots\\
&\ddots&\ddots
\end{pmatrix}\ =:\ A_{> m}\,,
\]
we have
\[
\spess A\ =\ \spess A_{> m}\ \subseteq\ \spec A_{> m}\ \subseteq\ \bigcup_{k\ge m}\Specn_{\eps_{n}} (AP_{n,k},A^*P_{n,k}),\qquad m,n\in\N,
\]
using the semi-infinite version of \eqref{eq:sandwich0} in the last step.
Taking the intersection over all $m,n\in\N$ gives
\begin{equation} \label{eq:spess1}
\spess A\quad\subseteq\quad
\bigcap_{n\in\N}\bigcap_{m\in\N}\bigcup_{k\ge m}\Specn_{\eps_{n}} (AP_{n,k},A^*P_{n,k}).
\end{equation}

In \cite{SCW-ML:Inclusion2}, using results from \cite{HagLiSei}, we prove the following:
\begin{proposition}{\bf (semi-infinite)} For bounded band operators $A$ on $\ell^2(\N)$, %one has
formula \eqref{eq:spess1} holds in fact with ``$\subseteq$'' replaced by equality. In addition, after this replacement,
\begin{itemize}\itemsep-1mm
%\item[\bf a) ] the subset sign  in \eqref{eq:spess1} is an equality;
\item[\bf a) ] the intersection sign ``$\cap_{n\in\N}$'' in \eqref{eq:spess1} can be replaced by a Hausdorff-limit $\lim_{n\to\infty}$;
\item[\bf b) ] the two intersection signs ``$\cap_{n\in\N}\cap_{m\in\N}$'' in \eqref{eq:spess1} can be replaced by a single Hausdorff-limit $\lim_{m=n\to\infty}$. 
\end{itemize}
\end{proposition}

%\begin{acknowledgement}
%  An acknowledgement may be placed at the end of the~article.
%\end{acknowledgement}

%\textcolor{red}{They actually want us to drop the article's title in journal papers, just giving authors, journal, volume, pages, year. But then we're in big trouble with \cite{SCW-ML:Inclusion1,SCW-ML:Inclusion2} -- so I put the title everywhere. Alternatively: put a version 1.0 of \cite{SCW-ML:Inclusion2} on arxiv.org now and refer to that?}

%\textcolor{red}{Remember: We might want to change the title of \cite{SCW-ML:Inclusion2}.}


\begin{thebibliography}{99}
\bibitem{CollComp1}
P.~M.~Anselone,
Collectively Compact Operator Approximation Theory and Applications to Integral Equations,
Prentice-Hall, Englewood Cliffs 1971.


\bibitem{AvronMoucheSimon}
J.~Avron, P.H.M.~v.Mouche and B.~Simon,
%On the measure of the spectrum for the almost Mathieu operator,
Communications in Mathematical Physics {\bf 132}, 103--118 (1990).

\bibitem{TriRand}
S.N.~Chandler-Wilde, R.~Chonchaiya and M.~Lindner,
%On the Spectra and Pseudospectra of a Class of Non-Self-Adjoint Random Matrices and Operators,
Operators and Matrices {\bf 7}, 739--775 (2013).

\bibitem{SCW-ML:Inclusion1}
S.N.~Chandler-Wilde, R.~Chonchaiya and M.~Lindner,
On Spectral Inclusion Sets and Computing the Spectra and Pseudospectra of Bounded Linear Operators,
in preparation.

\bibitem{SCW-ML:Inclusion2}
S.N.~Chandler-Wilde and M.~Lindner,
Convergent spectral inclusions for finite and infinite banded matrices,
in preparation.

\bibitem{CollComp2}
S.~N.~Chandler-Wilde and B.~Zhang,
%A generalised collectively compact operator theory with an application to second kind integral equations on unbounded domains,
J. Integral Equations Appl. {\bf 14}, 11--52 (2002).

\bibitem{HengPhD}
R.~Chonchaiya,
Computing the Spectra and Pseudospectra of Non-Self-Adjoint Random Operators Arising in Mathematical Physics, 
PhD thesis, Reading, 2010.

\bibitem{Elliot}
M.-D.~Choi, G.A.~Elliot and N.~Yui,
%Gauss polynomials and the rotation algebra,
Inventiones mathematicae {\bf 99}, 225--246 (1990).

\bibitem{Davies}
E.B.~Davies,
Linear Operators and their Spectra, 
Cambridge University Press, 2007.

\bibitem{subwords}
F.~Gabel, D.~Gallaun, J.~Großmann, M.~Lindner and R.~Ukena,
%Spectral approximation of generalized Schrödinger operators via approximation of subwords,
arXiv:2209.11613.

\bibitem{HagLiSei}
R.~Hagger, M.~Lindner and M.~Seidel,
%Essential pseudospectra and essential norms of band-dominated operators,
Journal of Mathematical Analysis and Applications {\bf 437}, 255--291 (2016).

\bibitem{LenzPhD}
D.~Lenz,
Aperiodische Ordnung und gleichmäßige spektrale Eigenschaften von Quasikristallen, 
PhD thesis, Logos Verlag Berlin 2000.

\bibitem{MLBook}
M.~Lindner,
Infinite Matrices and their Finite Sections: An Introduction to the Limit Operator Method,
Frontiers in Mathematics, Birkh\"auser 2006.

\bibitem{BigQuest}
M.~Lindner, M.~Seidel,
%An Affirmative Answer to a Core Issue on Limit Operators,
Journal of Functional Analysis {\bf 267}, 901--917 (2014).

\bibitem{RaRoSiBook}
V.S.~Rabinovich, S.~Roch and B.~Silbermann,
Limit Operators and Their Applications in Operator Theory,
Birkh\"auser 2004.

\bibitem{Yu}
H.~Sako,
%Property A and the operator norm localization property for discrete metric spaces,
Journal für die reine und angewandte Mathematik (Crelle’s Journal) {\bf 690}, 207--216 (2014).

\bibitem{Shargo}
E.~Shargorodsky,
%On the level sets of the resolvent norm of a linear operator,
Bull. London Math. Soc. {\bf 40}, 493--504 (2008).

\bibitem{TrefEmb}
L.N.~Trefethen and M.~Embree,
Spectra and Pseudospectra: The Behavior of Nonnormal Matrices and Operators,
Princeton University Press, Princeton, NJ, 2005.

\end{thebibliography}
\end{document}